\input amstex
\documentstyle{amsppt}
\tolerance 3000 \pagewidth{5.5in} \vsize7.0in
\magnification=\magstep1 \NoBlackBoxes \NoRunningHeads
\widestnumber \key{AAAAAAAAAAA} \topmatter
\author L. Brandolini, S. Hofmann, and A. Iosevich
\endauthor
\title Sharp rate of average decay of the Fourier transform of a
bounded set
\endtitle
\address L. Brandolini, Dipartimento di Ingegneria,
Universit\`a di Bergamo, V.le G. Marconi 4, 24044 Dalmine BG,
Italy \ email:brandolini \@unibg.it
\endaddress

\address S. Hofmann, Department of Mathematica, University of
Missouri, Columbia Missouri 65211, USA \ email: hofmann \@
math.missouri.edu \endaddress

\address A. Iosevich, Department of Mathematica, University of
Missouri, Columbia Missouri 65211, USA \ email: iosevich \@
wolff.math.missouri.edu
\endaddress

\date February 10, 2003
\enddate
\thanks Research supported in part by the NSF grants and INdAM
\endthanks

\abstract Estimates for the decay of Fourier transforms of
measures have extensive applications in numerous problems in
harmonic analysis and convexity including the distribution of
lattice points in convex domains, irregularities of distribution,
generalized Radon transforms and others. Here we prove that the
spherical $L^2$-average decay rate of the Fourier transform of the
Lebesgue measure on an arbitrary bounded convex set in ${\Bbb
R}^d$ is
$$
{\left(\int_{S^{d-1}} {|\widehat{\chi}_B(R\omega)|}^2 d\omega
\right)}^{\frac{1}{2}} \lesssim R^{-\frac{d+1}{2}}. \tag*
$$
This estimate is optimal for any convex body and in particular it
agrees with the familiar estimate for the ball. The above estimate
was proved in two dimensions by Podkorytov, and in all dimensions
by Varchenko under additional smoothness assumptions. The main
result of this paper proves $(*)$ in all dimensions under the
convexity hypothesis alone. We also prove that the same result
holds if the boundary of $\partial \Omega$ is $C^{\frac{3}{2}}$.
\endabstract

\endtopmatter
\document
\newpage

\head Introduction \endhead

\vskip.125in

Let $B$ be a bounded open set in ${\Bbb R}^d$. If $\partial B$ is
sufficiently smooth and has everywhere non-vanishing Gaussian
curvature, then
$$ |\widehat{\chi}_B(R \omega)| \lesssim R^{-\frac{d+1}{2}},
\tag0.1$$ with constants independent of $\omega$, where
$$ \hat{f}(\xi)=\int_{{\Bbb R}^d} e^{-i x \cdot \xi}
f(x)dx, $$ denotes the Fourier transform, and $A \lesssim B$ means
that there exists a positive constant $C$ such that $|A| \leq
C|B|$. The estimate $(0.1)$ is optimal in a very strong sense. One
can check that a better rate of decay at infinity is not possible.
One can also check that if the Gaussian curvature vanishes at even
a single point, then $(0.1)$ does not hold.

In fact, the point-wise estimate may be much worse. For example,
if $B$ is convex, one has
$$ |\widehat{\chi}_B(R \omega)| \lesssim R^{-1}, $$ and
the case of a cube ${[0,1]}^d$ shows that one cannot, in general,
do any better. See, for example, \cite{St93}, for a nice
description of these classical results.

In spite of the fact that the estimate $(0.1)$ does not hold in
general, a basic question is whether this estimate holds on
average for a large class of domains, for example, bounded open
sets with a rectifiable boundary. More precisely, one should like
to know for which domains one has the following estimate:
$$ {\left( \int_{S^{d-1}} {|\widehat{\chi}_B(R \omega)|}^2
d\omega \right)}^{\frac{1}{2}} \lesssim R^{-\frac{d+1}{2}}.
\tag0.2$$

In some cases, it is equally useful to know whether
$$ {\left( \int_{S^{d-1}} {|\widehat{\sigma}(R \omega)|}^2
d\omega \right)}^{\frac{1}{2}} \lesssim R^{-\frac{d-1}{2}},
\tag0.3$$ where $\sigma$ is the Lebesgue measure on the boundary
of $B$. Under a variety of assumptions, for example, if $\partial
B$ is Lipschitz, $(0.2)$ and $(0.3)$ are linked via the divergence
theorem. We use this fact in the proof of our main result below.

An example due to Sj\"olin (\cite{Sj93}) shows that $(0.3)$ is not
purely dimensional. He showed that if $\sigma$ is an arbitrary
$(d-1)$-dimensional compactly supported measure, then the best
exponent one can expect on the right hand side of $(0.3)$ is
$\frac{d-\frac{3}{2}}{2}$. This means that in order to prove an
estimate like $(0.2)$ we must use the fact that $\partial B$ is in
some sense a hyper-surface.

Several results of this type have been proved over the years. In
\cite{Pod91}, Podkorytov proved $(0.2)$ for convex domains in two
dimensions using a beautiful geometric argument that relied on the
fact that in two dimensions, the Fourier transform of a
characteristic function of a convex set in a given direction is
bounded by a measure of a certain geometric cap. See, for example,
\cite{BrNaWa88} or \cite{BRT98} for more details. Unfortunately,
in higher dimensions one cannot bound the Fourier transform of a
characteristic function of a convex set by such a geometric
quantity. See, for example, \cite{BMVW88}. For the case of average
decay on manifolds of co-dimension greater than one, see e.g.
\cite{Christ85}, \cite{Marshall88}, and \cite{IoSa97}.

The analytic case has been known for a long time. See, for
example, \cite{R66}. In \cite{Var83}, Varchenko proved $(0.2)$
under the assumption that $\partial B$ is sufficiently smooth.
Smoothness allows one to use the method of stationary phase in a
very direct and strong way. In the general case, one must come to
grips with the underlying geometry of the problem. In the main
result of this paper, we drop the smoothness assumption and prove
that $(0.2)$ holds for all bounded open convex sets $B$ in ${\Bbb
R}^d$. In addition, we prove the same estimate under an assumption
that the boundary is $C^{\frac{3}{2}}$.

The main geometric feature of our approach is a quantitative
exploitation of the following simple idea: if $\omega \in S^{d-1}$
is normal to $\partial B$ at $x$, and $y$ is sufficiently close to
$x$, then $x-y$ cannot be parallel to $\omega$. This allows us to
deal with the so-called "stationary" points of the oscillatory
integral resulting from $(0.3)$. Unlike the smooth case, where
"non-stationary" points are very easy to handle using integration
by parts, in the general case one is forced to exploit the
smoothness of the sphere along with an appropriate integration by
parts argument that exploits either convexity or the
$C^{\frac{3}{2}}$ assumption on the boundary.

The estimates $(0.2)$ and $(0.3)$ have numerous applications in
various problems of harmonic analysis, analytic number theory and
geometric measure theory. Moreover, $(0.2)$ and $(0.3)$ imply
immediate generalizations of a number of results in analysis and
analytic number theory to higher dimensions. See, for example,
\cite{BC00}, \cite{BCT97}, \cite{BRT98}, \cite{CdV77},
\cite{Christ85}, \cite{Hu96}, \cite{IoSa97}, \cite{KolWolff02},
\cite{Mat87}, \cite{Mont94}, \cite{R66}, \cite{RT01}, \cite{Sj93},
\cite{Sk98}, and \cite{Var83}. We give two simple examples to
illustrate the point.

\subhead Distribution of lattice points in convex domains
\endsubhead A classical result due to Landau says that if $B$ is
convex, and $\partial B$ is smooth and has non-vanishing Gaussian
curvature, then
$$ \left| \# \{tB \cap {\Bbb Z}^d\}-t^d|B| \right| \leq
Ct^{d-2+\frac{2}{d+1}}. \tag0.4$$

The proof is based on $(0.1)$. Using $(0.2)$ instead, one can
prove the following version of $(0.4)$:
$$ {\left( \int_{S^{d-1}} {\left| \# \{t\rho B \cap {\Bbb Z}^d \}-t^d|B|
\right|}^2 d\rho \right)}^{\frac{1}{2}} \leq
Ct^{d-2+\frac{2}{d+1}}, \tag0.5$$ where $\rho B$ denotes the
rotation of $B$ by $\rho \in S^{d-1}$ viewed as an element of
$SO(d)$. See, for example, \cite{Ios2001}, \cite{BCIPT02}, and
\cite{BIT01} for a detailed discussion of applications of average
decay of the Fourier transform to lattice point problems. Also
note that Theorem 1.2 below shows that convexity may be replaced
by a $C^{\frac{3}{2}}$ assumption.

\subhead Falconer Distance Problem \endsubhead A result due to
Falconer (\cite{Falc86}) says that if the Hausdorff dimension of a
set $E \subset {[0,1]}^d$, $d>1$, is greater than $\frac{d+1}{2}$,
then the distance set $\Delta(E)=\{|x-y|: x,y \in E\}$ has
positive Lebesgue measure.

Let $\Delta_B(E)=\{{||x-y||}_B: x,y \in E\}$, where
${||\cdot||}_B$ denotes the distance induced by a bounded convex
set $B$. Using $(0.2)$ one can prove that if the Hausdorff
dimension of a set $E \subset {[0,1]}^d$, $d>1$, is greater than
$\frac{d+1}{2}$, then $\Delta_B(\rho B)$ has positive Lebesgue
measure for almost every $\rho \in S^{d-1}$ viewed as an element
of $SO(d)$. One can also apply this technology to geometric
combinatorial results in a discrete setting. See, for example,
\cite{HoIo2002} and \cite{IoLa2002} for a detailed discussion of
these issues.

\head Section 1: $L^2$-Average decay  \endhead

\vskip.125in

Our main results are the following two theorems.

\proclaim{Theorem 1.1} Let $B$ be a bounded convex domain in
${\Bbb R}^d$. Then
$$ \int_{S^{d-1}} {|\widehat{\chi}_B(R \omega)|}^2 d\omega \lesssim
R^{-(d+1)}. \tag1.1 $$
\endproclaim

\proclaim{Theorem 1.2} Let $B$ be an bounded open set in ${\Bbb
R}^d$ satisfying the following assumption. The boundary of $B$ can
be decomposed into finitely many neighborhoods such that given any
pair of points $P,Q$ in the neighborhood,
$$ \left|(P-Q) \cdot n(Q) \right| \lesssim
{|P-Q|}^{\frac{3}{2}}, \tag1.2$$ where $n(Q)$ denotes the unit
normal to $\partial B$ at $Q$. Then $(1.1)$ holds.
\endproclaim

\head Proof of Theorem 1.1 and Theorem 1.2 \endhead

\vskip.125in

We shall give simultaneous proofs of Theorem 1.1 and Theorem 1.2.
The argument is based on the fact that both convex surfaces and
$C^{\frac{3}{2}}$ surfaces satisfy the following geometric
condition:

The boundary of $B$ can be decomposed into finitely many
neighborhoods $B_j$ such that on each neighborhood the surface is
given as a graph of a Lipschitz function with the Lipschitz
constant $<1$.

The geometric meaning of this condition is that for $x$ and $y$
belonging to the same neighborhood $B_j$, the secant vector $x-y$
lies strictly within $\frac{\pi}{4}$ of our local coordinate
system's horizon. We shall henceforth refer to this as the
``secant property".

This geometric condition is clearly satisfied by $C^1$ (and hence
$C^{\frac{3}{2}}$) surfaces, by taking the neighborhoods to be
sufficiently small. For convex domains, we make the following
construction. We cover $S^{d-1}$ by a smooth partition of unity
$\eta_j$, such that $\sum_j \eta_j \equiv 1$ and such that support
of each $\eta_j$ is contained in the intersection of the sphere
and a cone of aperture strictly smaller than $\frac{\pi}{2}$. Then
if $n(x)$ denotes the Gauss map taking $x \in \partial B$ to the
unit normal at $x$, then $\sum_j \eta_j(n(x))$ induces the desired
decomposition on the boundary of $B$. We note that in the convex
case the number of such neighborhoods depends only on dimension.

By the divergence theorem,
$$ \widehat{\chi}_B(R \omega)=-\frac{1}{2\pi iR}
\int_{\partial B} e^{-i x \cdot R \omega} \left( \omega \cdot n(x)
\right) d\sigma(x), \tag1.3$$ where $n(x)$ denotes the unit normal
to $\partial B$ at $x$, and $d\sigma$ denotes the surface measure
on the boundary. This reduces the problem to the boundary of $B$.

\subhead Decomposition of the boundary \endsubhead Let $\phi_j$
denote a smooth partition of unity on $\partial B$ subordinate to
the decomposition $\partial B=\cup_{j=1}^N B_j$. Moreover,
$\phi_j$s are chosen such that on the support of each $\phi_j$,
the aforementioned secant property still holds. It follows that
the corresponding Lipschitz constant $K_j$ is less than $1$. This
is the basic building block of our proof.

Let $\psi_j$ be a smooth cutoff function identically equal to $1$
on the spherical cap of solid angle $>\frac{\pi}{2}$ and which is
supported in a slightly bigger spherical cap which lies at a
strict positive distance from all the vectors $\frac{x-y}{|x-y|}$,
$x,y \in supp(\phi_j) \subset \partial B$. Notice that our
hypothesis make such a decomposition possible and that all the
vectors normal to $\partial B$ on the support of $\phi_j$ lie
strictly inside the support of $\psi_j$.

\subhead Singular directions \endsubhead This part of the proof is
identical in the convex and the $C^{\frac{3}{2}}$ cases. In fact,
it depends only on the secant property. Let
$$ F_j(R \omega)=\int_{\partial B} e^{-i x \cdot R \omega}
d\mu_j(x), $$ where
$$ d\mu_j=\left( \omega \cdot n(x) \right) \phi_j(x) d\sigma(x).
$$

In view of $(1.3)$ and the triangle inequality, it suffices to
show that
$$ \int_{S^{d-1}} {|F_j(R \omega)|}^2 d\omega \lesssim
R^{-(d-1)}.$$

Now,
$$ \int_{S^{d-1}} {|F_j(R \omega)|}^2 d\omega
=\int_{S^{d-1}} {|F_j(R \omega)|}^2 \psi_j(\omega)
d\omega+\int_{S^{d-1}} {|F_j(R \omega)|}^2
(1-\psi_j(\omega))d\omega=I+II.$$

We shall refer to the support of $\psi_j$ as ``singular"
directions, and the other vectors on the sphere as ``non-singular"
directions. The origin of this notation is the fact that in the
smooth case, the singular, or stationary directions are the ones
that are normal to the relevant piece of the hyper-surface in
question.

We have
$$ I=\int_{\partial B} \int_{\partial B} \int_{S^{d-1}} e^{i (x-y) \cdot R \omega} \psi_j(\omega) d\mu_j(x) d\mu_j(y).
$$

Using the definition of $\psi_j$, we integrate by parts $N$ times
and obtain
$$ I \lesssim \int_{\partial B} \int_{\partial B} \min \{1,
{(R|x-y|)}^{-N} \} d\mu_j(x) d\mu_j(y) \lesssim R^{-(d-1)},
$$ since $d\mu_j$ is $d-1$-dimensional and compactly
supported.

\subhead Non-singular directions \endsubhead We shall take the
following perspective on the spherical coordinates. Let
$\omega=\omega(\tau_1, \dots, \tau_{d-2}, \theta)$, where
$(\tau_1, \dots, \tau_{d-2})$ denotes the "azimuthal" angles, and
$\theta$ denotes the remaining angle, i.e
$\theta=\tan^{-1}(x_d/x_1)$. Note that for each fixed $\theta$,
$(\tau_1, \dots, \tau_{d-2})$ give a coordinate system on the
``great circle" tilted at the angle $\theta$ from the horizontal.

For each fixed $\theta$, we set up a coordinate system such that
$$ II=\int_{-\frac{\pi}{2}}^{\frac{\pi}{2}}
\int {\left| \int e^{iR \omega' \cdot u} \Phi_{\theta}(u) du
\right|}^2 J(\tau, \theta) (1-\psi_j(\omega)) d\tau d\theta, \tag
2.1$$ where $\omega=\omega(\tau, \theta)$, $\omega=(\omega',
\omega_d)$,
$$ \Phi_{\theta}(u,\omega)=\omega\cdot(\nabla A_\theta(u),-1)\ \phi_j(u, A_{\theta}(u)),
$$
and $J$ is the (smooth) Jacobian corresponding to the spherical
coordinates. Here we are viewing this portion of the boundary of
$B$ as the graph, of the function $A_{\theta}$, above the
hyperplane determined by the $(d-2)$-dimensional ``great circle"
obtained by fixing $\theta$. Observe that
$\Phi_{\theta}(u,\omega)$ is linear in $\omega$. This fact and the
Minkowski inequality allows to assume in the following $
\Phi_{\theta}(u,\omega)$ independent of $\omega$.

By a further partition of unity, a rotation, and the triangle
inequality, we may assume that we are in an arbitrarily small
neighborhood of $\omega=(1,0, \dots, 0)$.

The key object in the remaining part of the proof is the
difference operator
$$ \Delta_hf(s)=f(s+h)-f(s). $$

We observe that the transpose of this operator
$$ \Delta^{*}_h=\Delta_{-h}. $$

We also note that
$$ \Delta_{\frac{1}{R}}(e^{iR \omega_1 u_1})=(e^{i
\omega_1}-1)e^{iR \omega_1 u_1}.$$

Then by discrete integration by parts, the square root of the
portion of $(2.1)$ in the neighborhood of $(1,0, \dots, 0)$ equals
$$ {\left( \int_{-\frac{\pi}{2}}^{\frac{\pi}{2}}
\int {\left| \frac{1}{e^{i\omega_1}-1} \int \int e^{iR \omega'
\cdot u} \Delta_{-\frac{1}{R}} \Phi_{\theta}(\cdot,u')(u_1)
du_1du' \right|}^2 J(\tau, \theta) \Psi_j(\omega) d\tau d\theta
\right)}^{\frac{1}{2}}, \tag 2.2$$ where $\Psi_j$ is an
appropriate cut-off function supported in the neighborhood of
$(1,0, \dots, 0)$, and $u'=(u_2, \dots, u_{d-1})$.

Applying the Minkowski integral inequality, we see that $(2.2)$ is
bounded by
$$ \int {\left( \int_{-\frac{\pi}{2}}^{\frac{\pi}{2}}
\int {\left| \int e^{iR \omega'' \cdot u'} \Delta_{-\frac{1}{R}}
\Phi_{\theta}(\cdot,u')(u_1)du' \right|}^2 J(\tau, \theta)
\Psi_j(\omega) d\tau d\theta \right)}^{\frac{1}{2}} du_1, \tag
2.3$$ where $\omega''=(\omega_2, \dots, \omega_{d-1})$. For a
fixed $\theta$, the integration in $\tau$ is over the
$d-2$-dimensional ``great circle". We may parameterize the sphere
so that this ``great circle" is given by
$\omega_1=\omega_1(\omega'')$.

Expanding $(2.3)$ and rewriting, we get

\vskip.125in
$$ \int {\left( \int_{-\frac{\pi}{2}}^{\frac{\pi}{2}}
\int \int e^{iR \omega'' \cdot (u'-v')} \Delta_{-\frac{1}{R}}
\Phi_{\theta}(\cdot,u')(u_1)\Delta_{-\frac{1}{R}}
\Phi_{\theta}(\cdot,v')(u_1)du'dv' J'(\omega'', \theta)
\Psi_j(\omega) d\omega'' d\theta \right)}^{\frac{1}{2}} du_1, \tag
2.4$$ where $J'(\omega'', \theta)$ is smooth in $\omega''$.

Integrating by parts in $\omega''$ we see that $(2.4)$ is bounded
by
$$ \int {\left(
\int_{-\frac{\pi}{2}}^{\frac{\pi}{2}} \int \int \min \{1,
{R|u'-v'|}^{-N}\} |\Delta_{-\frac{1}{R}}
\Phi_{\theta}(\cdot,u')(u_1)||\Delta_{-\frac{1}{R}}
\Phi_{\theta}(\cdot,v')(u_1)|du'dv' q(\theta) d\theta
\right)}^{\frac{1}{2}} du_1 \tag 2.5$$
$$ \lesssim R^{-\frac{d-2}{2}} \int {\left(
\int_{-\frac{\pi}{2}}^{\frac{\pi}{2}} \int {|{\Cal M}'
\Delta_{-\frac{1}{R}} \Phi_{\theta}(\cdot,u')(u_1)|}^2 du'
q(\theta) d\theta \right)}^{\frac{1}{2}} du_1, $$ where ${\Cal
M}'$ is the Hardy-Littlewood maximal function in the $u'$
variable, and $q$ is a smooth cutoff function. The last inequality
uses the standard fact that convolution with a radial, integrable
and decreasing kernel is dominated by the maximal function. See,
for example, \cite{St70}, Chapter 3. Since our integrand is
compactly supported in the $u_1$ variable, we may apply
Cauchy-Schwarz to obtain that $(2.5)$ is bounded by
$$
R^{-\frac{d-2}{2}} {\left( \int_{-\frac{\pi}{2}}^{\frac{\pi}{2}}
\int \int {|\Delta_{-\frac{1}{R}} \Phi_{\theta}(\cdot,u')(u_1)|}^2
du' du_1 q(\theta) d\theta \right)}^{\frac{1}{2}},
$$ since the
Hardy-Littlewood maximal function is bounded on $L^2$.

The conclusion of the theorem now follows from the estimate
$$ {||\Delta_{\frac{1}{R}} \Phi_{\theta}||}_{L^2(du)} \leq
CR^{-\frac{1}{2}}. \tag2.6$$

Clearly $(2.6)$ holds if $\partial B \in C^{\frac{3}{2}}$, for in
that case $\Phi_{\theta} \in C^{\frac{1}{2}}$ with compact
support. In the convex case we interpolate between the estimates
$$ {||\Delta_{\frac{1}{R}}\Phi_{\theta}||}_{L^{\infty}(du)} \leq
C, \tag2.7$$ and
$$ {||\Delta_{\frac{1}{R}}\Phi_{\theta}||}_{L^1(du)} \leq CR^{-1},
\tag2.8$$ where $(2.7)$ holds because convex surfaces are
Lipschitz (hence $\Phi_{\theta}$ is bounded), and $(2.8)$ holds by
the mean-value theorem, Fubini theorem, and Gauss-Bonnet theorem
(for cross-sections) in the $u_1$ variable.

\vskip.125in \vskip.125in

\newpage

\Refs

\ref \key BMVW88 \by J. Bak, D. McMichael, J. Vance, and S.
Wainger \paper Fourier transforms of surface area measure on
convex surfaces in ${\Bbb R}^3$ \jour Amer. J. of Math. \vol 111
\yr 1988 \pages 633-668
\endref

\ref \key BC00 \by L. Brandolini, L. Colzani, \paper Decay of
Fourier transforms and summability of eigenfunction expansions
\jour Ann. Scuola Norm. Sup. Pisa \vol 29 \yr 2000 \pages 611-638
\endref

\ref \key BCT97  \by L. Brandolini, L. Colzani, and G. Travaglini
\paper Average decay of Fourier transforms and integer points in
polyhedra \jour Ark. Mat. \vol 35 \yr 1997 \pages 253-275
\endref

\ref \key BCIPT02 \by L.Brandolini, L.Colzani, A. Iosevich,
A.N.Podkorytov, G.Travaglini \paper Geometry of the Gauss map and
lattice points in convex domain \jour Matematika (to appear)
\endref

\ref \key BIT01 \by L. Brandolini, A. Iosevich, and G. Travaglini
\paper Average decay of the Fourier transform and geometry of the
Gauss map \jour preprint \endref

\ref \key BRT98 \by L. Brandolini, M. Rigoli, and G. Travaglini
\paper Average decay of Fourier transforms and geometry of convex
sets \jour Revista Mat. Iber. \yr 1998 \vol 14 \pages 519--560
\endref

\ref \key BrNaWa88 \by J. Bruna, A. Nagel and S. Wainger \paper
Convex hyper-surfaces and Fourier transforms \jour Ann. of Math.
\vol 127 \yr 1988 \pages 333-365 \endref

\ref \key CdV77 \by Y. Colin de Verdiere \paper Nombre de points
entiers dans une famille homothétique de domains de ${R}$ \jour
Ann. Sci. École Norm. Sup.  \vol 10 \yr 1977 \pages 559--575
\endref

\ref \key Christ85 \by M. Christ \paper On the restriction of the
Fourier transform to curves: endpoint results and the degenerate
case \jour Trans. Amer. Math. Soc. \vol 287 \yr 1985 \pages
223-238 \endref

\ref \key Falc86 \by K. J. Falconer \paper On the Hausdorff
dimensions of distance sets \jour Mathematika \vol 32 \pages
206-212 \yr 1986 \endref

\ref \key He62 \by C.S. Herz \paper On the number of lattice
points in a convex set \jour Amer. J. Math. \vol 84  \yr 1962
\pages 126--133
\endref

\ref \key Hl50 \by E. Hlawka \paper \"{U}ber Integrale auf
konvexen K\"{o}rpen, I \& II \jour Monatsh. Math. \vol 54 \yr 1950
\pages 1--36, 81--99
\endref

\ref \key HoIo2002 \by S. Hofmann and A. Iosevich \paper Falconer
conjecture in the plane for random metrics \jour Israel Journal of
Math. (submitted) \yr 2002 \endref

\ref \key Hu96 \by M. N. Huxley \paper Area, lattice points, and
exponential sums \jour Clarendon Press, Oxford \yr 1996 \endref

\ref \key Ios2001 \by A. Iosevich \paper Lattice points and
generalized diophantine conditions \jour Journal of Number Theory
\vol 90, no.1 \yr 2001 \pages 19-30 \endref

\ref \key IoLo2002 \by A. Iosevich and I. Laba \paper $K$-distance
sets and the Falconer conjecture \yr 2002 \jour (in preparation)
\endref

\ref \key IoSa97 \by A. Iosevich and E. Sawyer \paper Averages
over surfaces \jour Adv. in Math. \vol 132 \yr 1997 \pages 46-119
\endref

\ref \key KolWolff02 \by M. Kolountzakis and T. Wolff \paper On
the Steinhaus tiling problem \jour Mathematika \vol 46 \yr 2002
\pages 253-280 \endref

\ref \key Marshall88 \by B. Marshall \paper Decay rates of Fourier
transforms of curves \jour Trans. of AMS. \vol 310 \yr 1988
\endref

\ref \key Mat87 \by P. Mattila \paper Spherical averages of
Fourier transforms of measures with finite energy; dimension of
intersections and distance sets \yr 1987 \vol 34 \pages 207-228
\endref

\ref \key Mont94 \by H. Montgomery \paper Ten lessons on the
interface between analytic number theory and harmonic analysis \yr
1994 \jour CBMS Regional Conference Series in Mathematics,
American Mathematical Society \vol 84 \endref

\ref \key Pod91 \by A.N. Podkorytov \paper The asymptotic of a
Fourier transform on a convex curve \jour Vestn. Leningr. Univ.
Mat. \vol 24 \yr 1991 \pages 57--65 \endref

\ref \key R66 \by B. Randol \paper A lattice-point problem II
\jour Trans. Amer. Math. Soc. \vol 125 \yr 1966 \pages 101--113
\endref

\ref \key RT01 \by F. Ricci and G. Travaglini \paper Convex
curves, Radon transforms and convolution operators defined by
singular measures \jour Proc. Amer. Math. Soc. \vol 129 \yr 2001
\pages 1739-1744 \endref

\ref \key Sj93 \by P. Sjolin \paper Estimates of spherical
averages of Fourier transforms and dimensions of sets \jour
Mathematika \yr 1993 \vol 40 \pages 322-330 \endref

\ref \key So93 \by C. Sogge \paper Fourier integrals in classical
analysis \jour Cambridge Univ. Press \vol 105 \yr 1993
\endref

\ref \key Sk98 \by M. M. Skriganov \paper Ergodic theory on
$SL(n)$, Diophantine approximations and anomalies in the lattice
point problem \jour Invent. Math. \vol 132 \yr 1998 \pages 1--72
\endref

\ref \key St70 \by E. M. Stein \paper Singular integrals and
differentiability properties of functions \jour Princeton
Mathematical Series, No. 30 Princeton University Press, Princeton,
NJ \yr 1970 \endref

\ref \key St93 \by E. M. Stein \paper Harmonic Analysis \yr 1993
\jour Princeton University Press \endref

\ref \key Str89 \by R. Strichartz \paper Harmonic analysis as
spectral theory of laplacians \jour J. of Func. Anal. \vol 87 \yr
1989 \page 51-148 \endref

\ref \key Str91 \by R. Strichartz \paper Fourier asymptotics of
fractal measures \jour J. Funct. Anal. \vol 89 \yr 1990 \pages
154-187 \endref

\ref \key Var83 \by A. Varchenko \paper Number of lattice points
in families of homothetic domains in ${\Bbb R}^n$ \jour Funk. An.
\vol 17 \yr 1983 \pages 1-6 \endref

\ref \key WZ77 \by R. Wheeden and A. Zygmund \paper Measure and
Integral \jour Pure and Applied Mathematics \yr 1977 \endref

\ref \key Wolff00 \by T. Wolff \paper Math 191 lecture notes \jour
Harmonic analysis lecture notes; CalTech \yr 2000 \endref

\endRefs
\enddocument